\documentclass[18pt]{article}
\usepackage[cp1251]{inputenc}
\usepackage[T2A]{fontenc}
\usepackage[english]{babel}
\usepackage {amsfonts}
\begin{document}
\begin{center}
{\bf The Problem of Projecting the Origin of Euclidean Space onto the Convex Polyhedron}\\
\smallskip

{\bf Z.~R.~Gabidullina }\\
\end{center}

\begin{center}
Kazan Federal University

{\bf Email}: zgabid@mail.ru, Zulfiya.Gabidullina@kpfu.ru \\
\end{center}

\newtheorem{mydef}{Definition}[section]

\newtheorem{theorem}{Theorem}[section]
\newtheorem{seqv}{Corollary}[section]
\newtheorem{note}{Remark}[section]
\newtheorem{lemm}{Lemma}[section]
\newtheorem{examp}{Example}[section]
\newtheorem{abstr}{Abstract.}
\newenvironment{proof}
    {\par\noindent{\bf Proof.\enspace}}
    {\hfill$\square$}

{\bf Abstract}. This paper is aimed at presenting a systematic
survey of the existing now different formulations for the problem
of projection of the origin of the Euclidean space onto the convex
polyhedron (PPOCP). In the present paper, there are investigated
the reduction of the projection program to the problems of
quadratic programming, maxi\-min, li\-near complementarity, and
nonnegative least squares. Such reduction justifies the
opportunity of utilizing a much more broad spectrum of powerful
tools of mathe\-matical programming for solving PPOCP.
 The paper's goal is to draw the attention
of a wide range of research at the different formulations of the
projection  problem, which remain largely unknown due to the fact
that the papers (addressing the subject of concern) are published
even though on the adjacent, but other topics, or only in the
conference proceedings.

{\bf Keywords}: projection,  convex polyhedron,
 quadratic programming,  ma\-xi\-min problem,
complementarity problem,  nonnegative least squares problem \\

 MSC 2010: 90C30, 90C25, 90C20,
90C33, 65K05

\section{Introduction}
The problem of projecting the origin of the Euclidean space onto a
convex polyhedron plays an invaluable role in the differentiable
and nondifferentiable convex optimization,  in the theory of
linear separation the convex polyhedra, data classification and
identification. For many years, the numerous applications of the
 problem of projecting the origin onto a convex polyhedron
motivates intensive studying the possibility of reducing this
problem (in the original formulation with an abstract constraint)
to the other related ones of mathematical programming
\cite{G_Wolfe}--\cite{G_Gabid6}. Such a reduction, as a rule,
tends to encourage the appearance of new methods for solving the
original problem and finding new applications for the already
existing methods.

 Indeed, the class of projection techniques is
successfully used for solving a variety of the theoretical and
applied mathematical problems.
 The methods  applying the concept of projection
are commonly referred to as  the projection methods
\cite{G_Cegielski_Censor}. These algorithms use different kinds of
projections onto the convex sets (including orthogonal projections
which are the subject of this article), and serve for solving the
optimization problems as well as the  so-called feasibility
problems (the problems of finding the point belonging to some
set). The adequate for the most applications of the projection
methods up-to-date overview of the literature is presented in
 \cite{G_Cegielski_Censor}--\cite{G_Bauschke_Borwein}. In these
sources, there is annotated the bibliography on the iterative
projective techniques that use a projection onto the convex and
closed sets specified by a finite family of sets (for example, the
intersection of the individual sets from this family).

The computational success of these methods of projection  is
determined by the fact that, as a rule, it is always easier to
implement the projection onto the individual sets with a simpler
structure than, for example, at their intersection.

Finding the projection of points of the space onto a convex set,
including a convex polyhedron, is used in particular in the method
of cyclic projections. The cyclic orthogonal projection method can
serve, for example, to find a point belonging to the set
 $\Phi = \Phi_{1} \cap
\Phi_{2}\cap \ldots \cap \Phi_{N},$ where $\Phi_{i},$ \, $i =
\overline{1,N}$ are the polyhedra specified by the convex hulls of
the finitely many vectors
\cite{G_Bauschke_Borwein}--\cite{G_Bregman}. We further note that
the problem of projection of some point  $p \in \mathbb{R}^{n}$
onto the convex polyhedron $L$ may be reduced to the problem of
projecting the origin onto the Minkowski difference $L-p$,\, as it
was shown in \cite{G_Gabid2} (see p.565). In
\cite{G_Gabid}--\cite{G_Gabid6},
\cite{G_Gabid2}--\cite{G_Gabid7},\, there are investigated the
applications of the problem under study  for solving the related
problems such as
   \begin{itemize}
      \item [$-$] The problem of the strong linear separation of the convex polyhedra,
      \item [$-$] The problem of determining the distance between the convex polyhedra,
      \item [$-$] The problem of concretizing the nearest points of the convex polyhedra,
      \item [$-$] The problem of estimating the thickness of the separator between the two convex polyhedra,
      \item [$-$] The problem of solving the system of the linear (or quadratic) inequalities,
      \item [$-$] The variational inequalities associated with the linear separation of the convex polyhedra.
     \end{itemize}
      PPOCP is also widely used  in nonsmooth convex
optimization \cite{G_Shor} on the whole and the convex minimax
problems in particular \cite{G_Dem1}. In minimization of the
nondifferentiable convex functions, the learnt problem is used
  for the choice of the descent direction as well as for
verification the stop criterion of optimization algorithms. With
the same purpose, there can be utilized PPOCP in constrained
minimization of the differentiable pseudo-convex function (see,
for instance \cite{G_Gabid10}). Due to a well known connectedness
of PPOCP with the program of computing the distance between the
convex polyhedra (through the Minkowski difference), all
enumerated in \cite{G_Gilbert} applications such as computer
graphics, robotics, collision detection, and  path finding in
presence of obstacles are closely relevant to the problem under
consideration. The role of PPOCP is also valuable for the
mathematical methods in biomedical imaging and intensity-modulated
radiation therapy \cite{G_Censor_Jiang_Louis}.

In present paper, our goal is to provide a systematic exposition
of the existing now various settings for PPOCP. Here, we study the
questions related to the issue of reducing the projection program
to the problems of quadratic programming, maximin, linear
complementarity, and nonnegative least squares.
 Our prime aim is to draw the attention
of a broad range of research at the different formulations of the
projection  problem, which remain largely unknown due to the fact
that the papers (addressing the subject of concern) are published
even though on the relevant, but other topics, or only in the
conference proceedings.

      \section{Projection problem in Original Setting}

%\section*{\centerline {$1^{0}. \quad \mbox{}$}}

In this section of the paper,  we study the considered projection
problem in classical setting.

 Let the polyhedron is represented as convex hull of finitely many vectors from the $n-$dimensional Euclidean space: \\$L:= conv \bigl
\{z_{i}\bigr \}_{i \in
  I}$,\, where
  \,$I= \bigl \{1,2,\cdots, m \bigr \}$,\,
  , i.e.
$$L= \Bigl \{ z\in \mathbb{R}^{n}: z=\sum \limits_{i \in I} \alpha_{i}z_{i} \Bigm|
\,\alpha \in \mathbb{R}^{m}_{+},
 \, \sum \limits_{i \in I}\alpha_{i}=1 \Bigr \}.$$
 By construction, the polyhedron \,$L \subset \mathbb{R}^{n}$\, is convex and compact.

 Let us recall a widely known and traditional formulation of\, PPOCP:
 \begin{equation}\label{t_1}
    \min \limits_{z\in L} \,\, \Bigl( f(z)=\|z\|^{2} \Bigr) ,
\end{equation}
(see, for instance, ~\cite{G_Dem2}, \cite{G_Dem1}).

Owing to convexity and closedness of \,$L$,\, and strict convexity
of the objective function, a solution of the problem~(\ref{t_1})
exists, and moreover it is unique (see, for instance, Theorem 1
from \cite{G_Vas}, p. 193).

 If \,$\min \limits_{z\in L} \,\,
\|z\|^{2}=\|\rho\|^{2}$,\, then the point \,$\rho \in L$\, is used
to call a projection of the origin of Euclidean space
\,$\mathbb{R}^{n}$\, onto \,$L$,\, and  namely the
problem~(\ref{t_1}) is commonly referred to as a problem of
projecting the origin onto the convex polyhedron \,$L$.\,

An optimality criterion for the solution of the
problem~(\ref{t_1}) is performed in the following statement.
\begin{lemm}\label{A_2}
For \,$\rho \in L$\, to be a nearest to the origin of
\,$\mathbb{R}^{n}$\, point of the polyhedron \,$L \subset
\mathbb{R}^{n}$,\, it is necessary and sufficient to hold the
following variational inequality:
\begin{equation}\label{v_1}
\langle z-\rho ,\rho \rangle \geq 0 \quad \forall z \in L.
\end{equation}
\end{lemm}
\proof  Due to the known theorem specifying a necessary and
sufficient condition for the smooth convex function \,$f(z)$\, to
attain a minimum on the convex set \,$L \subset \mathbb{R}^{n}$\,
at the point $\rho$,\, it is fulfilled the following correlation:
\,$\langle f'(\rho),z-\rho \rangle \geq 0 \quad \forall z \in
L.$\, $\square$

The previous lemma demonstrates that the problem~(\ref{t_1}) can
be reduced to the problem of finding the point \,$\rho \in L$\,
satisfying to the variational inequality~(\ref{v_1}).
\begin{lemm}\label{A_3}
 For the point \,$y \in \mathbb{R}^{n}$\,
to satisfy the variational inequality
\\\centerline{ \,$\langle z-y,\,y \rangle \geq 0 \quad \forall z
\in L$,\,} it is necessary and sufficient for the vector \,$y$\,
to be a solution of the system of the quadratic inequalities:
\\\centerline{ \,$\langle z_{i}-y,\,y \rangle \geq 0 \quad
\forall i \in I.$\,}
\end{lemm}

The assertion of the \, Lemma \,\ref{A_3}\, can be easily proved
by taking into account that, on the one hand, any point
 \,$z = z_{i}, \, i = \overline{1,m}$\,\, belongs to \,$L.$\, On
 the other hand, any point
 \,$z \in L$\,  can be represented as a convex combination of the
 vectors \,$z_{i}, \, i=\overline{1,m}.$\,

 By help of Lemmas~\ref{A_2}--
 \ref{A_3},\, we can reformulate the optimality criterion for the
 solution of problem~(\ref{t_1}) in the following constructive (computationally realized) form: \\
The vector \,$\rho \in L$\, is a nearest to the origin of
\,$\mathbb{R}^{n}$\, point of \,$L$\, if and only if the following
system of the quadratic inequalities is consistent:
\begin{equation}\label{E_1}
   \langle z_{i}-\rho,\rho \rangle \geq
0 \quad \,\forall i \in I.\,
\end{equation}

We underline that the formulation of the problem (\ref{t_1}) (in
contrast to the other formulations of PPOCP) has an abstract
constraint \,$z \in L$\, which does not cause inconvenience to the
development of some optimization framework. In this regard, there
should be mentioned the  classical nearest point algorithms such
as Gilbert's method \cite{G_Gilbert},  MDM-method \cite{G_Dem2},
and more recent suitable affine subspace method for the case of a
simplex \cite{G_Nurmi}.

\section{Reduction to the Quadratic Programming Problems}
 This section is devoted to a treatment of how PPOCP can be posed as the different quadratic programming problems (QPP).
 For solving of QPP, there is a wide range of the effective algorithms \cite{G_Daugavet}--\cite{G_Gould_Toint}.
  Let us note that QPP may be also solved by using a package
Optimization Toolbox of MATLAB.   For instance, the function
quadprog allows to minimize the quadratic objective function
subject to the linear restrictions.

 The rest of this section is
organized as follows.  In subsections 3.1--3.2, we collect the QPP
having the obvious close association with the cone of generalized
strong support vectors \cite{G_Gabid6} of the given polyhedron
$L$.\, Subsection 3.3 deals with the optimization problem which
consists in finding the convex combination of the vectors $z_{i},
i = \overline{1,m}$\, with the least norm.

 \subsection{}In this subsection of the paper, we will justify a reduction of
the problem (\ref{t_1})  to the following quadratic programming
problems formulated in~\cite{G_Gabid}:
    \begin{equation} \label{max1}
\max \limits_{y \in \Omega} \,\|y\|^{2}.
    \end{equation}
     \begin{equation} \label{min1}
\min \limits_{y \in D} \,\|y\|^{2},
     \end{equation}
$$\mbox{where}\enspace D=\Bigl \{ y\in R^{n}: \langle z_{i},y \rangle \geq 1, \,
i=\overline{1,m} \Bigr\}, \,\Omega =\Bigl \{y\in R^{n}: \langle
z_{i},y \rangle \geq  \|y\|^{2},\, i =\overline{1,m} \Bigr \}.$$
We further illustrate of how one can solve PPOCP by means of
normalizing of the best strong support vectors from the sets $D$
and $\Omega$. Let us remark that unlike the polyhedron $L$ having
the inner description, the set $D$ is given by the outer
representation.
 First we also note that, by construction, the set \,$D$\, is
convex and closed. It is well known in optimization theory that a
solution for the problem of minimizing the strictly convex
function on the nonempty convex and closed set exists, and its
uniqueness is justified (see, for instance, ~\cite{G_Karman},\,p.
41).\, In~\cite{G_Gabid},\, there was proved that the set
\,$\Omega$\, is an intersection of the \,$m$\,\enspace
\,$n$\,-dimensional balls of radii
$$\hspace{-1cm} R_{i}= \frac{\|z_{i}\|}{2}\quad \mbox{
 with centers at the points}\quad O_{i}= \frac{z_{i}}{2}.$$ Consequently,
the set \enspace $\Omega$ \enspace is convex and compact. Note
that, by construction of the set \,$\Omega$,\, it always holds:
\,${\bf 0} \in \Omega$.\, According to the formula~(\ref{E_1}),\,
we obviously have \,$\rho \in \Omega$.\,

The following two statements establish a natural and useful
connection  between the feasible and optimal solutions
of~(\ref{max1})\, and (\ref{min1}).
\begin{lemm}\label{l3}
If \,$y^{*} \in\Omega$,\, \,$y^{*}\neq{\bf 0}$;\enspace $\bar{y}
\in D$,\,
$$\mbox{then:\, 1) \,} \enspace \hat{y}\in \Omega \mbox{\enspace with respect to }\enspace \hat{y}= \frac {\bar{y}}{\|
\bar{y}\|^{2}}, \enspace 2)\mbox{\enspace\,}\enspace \tilde{y}\in
D \mbox{\enspace with respect to }\enspace \tilde{y}= \frac
{y^{*}}{\|y^{*}\|^{2}}.$$
\end{lemm}

The lemma can be proved by help of substitution.
\begin{lemm}\label{l4}
If $y^{*}\neq {\bf 0}$  is a solution of the problem~(\ref{max1}),
\,$\bar{y}$\, is a solution of the problem~(\ref{min1}), then
$$\phantom{m}\mbox{1)\,the vectors\,}\enspace \hat{y}= \frac {\bar{y}}{\|
\bar{y}\|^{2}}\enspace \mbox{and}\enspace \tilde{y}= \frac
{y^{*}}{\|y^{*}\|^{2}}\enspace \mbox{are the solutions of the
problems~(\ref{max1})},\, (\ref{min1}),$$\phantom{m}respectively;
\\ \phantom{m}2) a solution of~(\ref{max1}) is unique.
\end{lemm}
\proof Using Lemma~\ref{l3}, one can prove the first two
statements by reductio ad absurdum. Taking into account the
uniqueness of the solution to the problem~(\ref{min1}), by help of
the first assertion of this lemma, it is easy to verify that the
solution of~(\ref{max1})\, is unique, too. \,$\square$\,

 Due to Lemmas~\ref{l3}--\ref{l4}, by reductio ad absurdum,  it is not hard to prove the following statement.
\begin{lemm}\label{l5}\, If the vector $y^{*}= {\bf 0}$\, is the solution to the problem~(\ref{max1}),\, then the problem~(\ref{min1})\, has no solutions.
If the domain of the feasible solutions for~(\ref{min1})\, is
empty, then $y^{*}= {\bf 0}$\, is the solution to the
problem~(\ref{max1}).\,
\end{lemm}
\begin{lemm}\label{l7}\,
If \,$y^{*}$\, is the solution for the problem~(\ref{max1}), then
\,$\rho=y^{*}.$\,
\end{lemm}
The proof of this assertion can be found in~\cite{G_Gabid},\,
p.24.
\begin{lemm}\label{l8}\,
\,$y^{*}={\bf 0}$\, is the solution to the problem~(\ref{max1}) if
and only if \,${\bf 0} \in L$.\,
\end{lemm}
\proof  Since it holds \,$\rho ={\bf 0}$\, if and only if т\,${\bf
0} \in L$;\, then the statement of the lemma follows  directly
from the Lemma~\ref{l7}. \,\,$\square$\,

 The assertions of
Lemmas~\ref{l5},\, \ref{l8}\, implies the following statement.
\begin{lemm}\label{l9}\,
The domain of the feasible solutions for the problem~(\ref{min1})
is empty if and only if \,${\bf 0} \in L$.\,
\end{lemm}

Therefore, for solving the problem of projecting the origin of the
space onto the convex polyhedron, it suffices to find the solution
\,$y^{*}$\,  for the problem~(\ref{min1})\, or \,$\bar{y}$\,
for~(\ref{max1}).\, The assertion of Lemma~\ref{l7},\, thereby
provides that \,$\rho=y^{*}$.\, Due to
Lemmas~\ref{l7},\,\ref{l4},\, the formula
\,\,$\rho=\bar{y}/\|\bar{y}\|^{2}$\, is made obvious.

\subsection{} In this subsection, there is studied the reduction of the
projection problem to the dual one for the problem~(\ref{min1}).
We rewrite first the constraints of the problem~(\ref{min1})\, in the vector-matrix form as follows:\\
\begin{equation}\label{Cye_1}
min \, \|y\|^{2}
\end{equation}
\begin{equation}\label{Cye_2}
Cy+e \leq 0,
\end{equation}
where $C$ is $m \times n$ matrix, the rows of which are the
vectors $-[z_{i}]^{T}, i \in I; e$ is  $m-$dimen\-sio\-nal vector,
$e^{T}=(1, \cdots, 1).$\, We further construct the dual program
to~(\ref{Cye_1})--(\ref{Cye_2}).

Toward this end, we construct the Lagrange function for the
problem~(\ref{Cye_1})--(\ref{Cye_2}):\,
$$\Lambda(y,u)= \|y \|^{2}+\langle u,Cy+e \rangle, y \in \mathbb{R}^{n}, u \in \mathbb{R}^{m}_{+}.$$
In what follows, we use the following notation:
 \,$h(u)=\min \limits_{y \in \mathbb{R}^{n}} \, \Lambda(y,u).$\,
Then, the program \,$\max \limits_{u \in
\mathbb{R}^{m}_{+}}\,\,h(u)$\, is called a dual program
to~(\ref{Cye_1})--(\ref{Cye_2}), and  we will call its variables
\,$u_{1}, \ldots , u_{l}$\, the dual ones.  The fact, that
$\,\Lambda(y,u)\,$ is a strictly convex and quadratic function
with respect to $y$, provides the attainment of its minimum at the
point $\,y(u)\,$ such that:
$$\Lambda'_{y} (y,u)= 2y+C^{T}u=0 \Rightarrow y(u)= -\frac{1}{2} C^{T}u.$$
Therefore, the function $\, \Lambda(y(u),u)\,$ will have the form:
$$ \Lambda(y(u),u)=-\frac{1}{4} \langle u, C C^{T}u\rangle  +\langle e,u \rangle.\,$$
 For this reason, the dual program may be expressed as follows (see, for instance, \cite{G_Gabid}):
\begin{equation} \label{lagr1}
\min \limits_{u \geq 0} \, (f(u)= \frac{1}{4} \langle u, C
C^{T}u\rangle  -\langle e,u \rangle).
\end{equation}
Let us note that this expression of the dual problem has an
interesting role in establishing the connection of the original
problem with the nonnegative least squares problem and linear
complementarity one, discussed below.
\begin{theorem} \label{prime_dual}
For the vector \,$y^{*} \in D$\, to be a solution of the
program~(\ref{Cye_1})--(\ref{Cye_2}), it is necessary and
sufficient to exist the vector \,$u^{*} \in \mathbb{R}^{m}_{+}$\,
satisfying
$$2y^{*}+C^{T}u^{*} ={\bf 0},$$
 $$u^{*}(Cy^{*}+e)= {\bf 0}.$$
\end{theorem}
The previous theorem represents the application of Theorem~10.4
from~\cite{G_Suhar} (p.180) for the concrete setting
of~(\ref{Cye_1})--(\ref{Cye_2}).
\begin{theorem}\label{unb_hu} For the system~(\ref{Cye_2}) to be inconsistent, it is necessary and sufficient
  that the function \,$h(u)$\, be unbounded from above.
\end{theorem}
\proof Since the domain of the feasible solutions for the
problem~(\ref{lagr1}) is nonempty, the statement of this theorem
immediately follows from the second statement of Theo\-rem~10.3
from ~\cite{G-Gavurin}, p.115. \,$\square$\,

Theorem~\ref{unb_hu} together with Lemma~\ref{l9}  yield that the
function \,$f(u)$\, is unbounded from below if and only if \,${\bf
0} \in L.$\,

For finding $\,y^{*},\,$  we need to solve the problem
(\ref{lagr1}), which consist of minimizing the quadratic convex
function on the set with a simple structure, or equivalently to
solve the following problem:
$$ \max \limits_{u \geq 0} \, \Lambda(y(u),u).\,$$
By help of the next theorem, it is easily be shown the convexity
of the objective function for the problem~(\ref{lagr1}).
\begin{theorem}\label{t8} (\cite{G_Vas}, p. 173)
Let  $\,U\,$ be a convex set in \,$\mathbb{R}^{n}$,\, \,$int(U)
\neq \emptyset$,\, \\$f(u) \in C^{2}(U)$.\, Then, for function
$\,f(u)\,$ to be convex on $\,U,\,$ it is necessary and sufficient
to hold
$$\langle f''(u) \xi, \xi \rangle \geq 0$$
for all $\xi =(\xi^{1}, \cdots , \xi^{n}) \, \in \mathbb{R}^{n}\,$
and $u \in U.$
\end{theorem}
\begin{lemm}\label{G_fu}
    The function \,$f(u)$\, is convex on the set \,$U=\mathbb{R}^{m}_{+}.$\,
\end{lemm}
\proof  The direct calculations with the following formulas:
$$f'(u)=\frac{1}{2}C C^{T}u -e, \enspace f''(u)=\frac{1}{2}C C^{T}$$
gives us
$$\langle f''(u) \xi, \xi \rangle = \frac{1}{2}\langle C C^{T} \xi, \xi \rangle= \frac{1}{2}\|C^{T} \xi\|^{2} \geq 0 \enspace
\mbox{\, for all\,} \, \xi\in \mathbb{R}^{m}, u \in
\,\mathbb{R}^{m}_{+}.$$ From Theorem~\ref{t8},\, it obviously
follows that the objective function \,$f(u)$\, of the
problem~(\ref{lagr1})\, is convex in \,$\mathbb{R}^{m}_{+}$.\,
$\square$

Let the problem~(\ref{lagr1}) be solved, and the vector
\,$u^{*}$\, be its solution. Then, according to
Theorem~\ref{prime_dual},\, one may determine the solution of the
prime problem~(\ref{Cye_1})--(\ref{Cye_2}) by the explicit formula
  \begin{equation}\label{y*}
\bar{y}= -\frac{1}{2} C^{T}u^{*}.
\end{equation}
It remains only to calculate the solution of the
problem~(\ref{t_1}) by the formula
\,$\rho=\bar{y}/\|\bar{y}\|^{2}.$\,

\subsection{}
In the subsection, a subject of our investigation is the
optimization problem often used for the construction of the data
classifier in machine learning \cite{G_Bennett_Mangasarian}. To
the best of our know\-ledge, the use of this optimization problem
was pioneered by Wolfe \cite{G_Wolfe} for finding the closest to
the origin of $\mathbb{R}^{n}$\, point on the polyhedron \,$L =
conv\{z_{i}\}_{i \in I}$. The projection problem  in this setting
becomes especially of interest later in the development of SVM
(Support Vector Machines). It is widely used for determining the
nearest points of the convex hulls for the two given sets of the
points (see, for instance,~\cite{G_Bennett}).

 We now consider the  projection problem in above-mentioned
setting:
\begin{equation} \label{G_norm1}
\min\,\, \Bigl( \varphi(\alpha)= \| \sum \limits_{i \in I}
\alpha_{i}z_{i}\|^{2} \Bigr),
\end{equation}
\begin{equation} \label{G_norm2}
\sum \limits_{i \in I} \alpha_{i} =1,
\end{equation}
\begin{equation} \label{G_norm3}
\alpha_{i}\geq 0,\,i \in I.
\end{equation}
Convert the problem~(\ref{G_norm1})--(\ref{G_norm3}) to a standard
vector-matrix form:
\begin{equation} \label{G_norm11}
\min\,\, \Bigl( \varphi(\alpha)= \langle \alpha, B \alpha \rangle
\Bigr),
\end{equation}
\begin{equation} \label{G_norm21}
\langle e, \alpha \rangle = 1,
\end{equation}
\begin{equation} \label{G_norm31}
\alpha \geq 0,
\end{equation}
where \,B\, is  \,$m \times m$--dimensional matrix,\, whose
element indexed by \,$j$\, in the row indexed by \,$i$\, is equal
to \,$\langle z_{i},z_{j} \rangle$,\,\, $i \in I$,\,\, $j \in
I,$\, $e^{T}= \underbrace{(1, \cdots, 1)}_{m}.$\, Obviously,
\,$B=C\cdot C^{T}.$\, This property of the matrix \,$B$\, allows
to reveal that the function \,$\varphi(\alpha)$\, is convex.
\begin{lemm}\label{A_12}
    The problem~(\ref{G_norm11})--(\ref{G_norm31}) is the convex quadratic programming problem.
\end{lemm}
\proof Since \,$\varphi''(\alpha)=2B,$\, it holds $$\langle
\varphi''(\alpha) \xi, \xi \rangle = 2\langle B \xi, \xi
\rangle=2\langle C C^{T} \xi, \xi \rangle= 2\|C^{T} \xi\|^{2}\geq
0 \enspace \forall \xi \in \mathbb{R}^{m}, \, \forall \alpha \in
\mathbb{R}^{m}.$$ Then, due to Theorem~\ref{t8},\, the function
\,$\varphi(\alpha)$\, is convex in \,$\mathbb{R}^{m}_{+}.$\,
Consequently, \,$\varphi(\alpha)$\, is also convex on the convex
set \,$\{\alpha \in \mathbb{R}^{m}_{+}: \langle e, \alpha \rangle
= 1 \}$.\, \,$\square$

By virtue of the non-emptiness of the set \,$\{\alpha \in
\mathbb{R}^{m}_{+}: \langle e, \alpha \rangle = 1 \}$\, and
boundedness from below of the function \,$\varphi(\alpha)$\, on
this set,
 the problem~(\ref{G_norm11})--(\ref{G_norm31})\,
 is solvable in the case of any convex polyhedron
\,$L =~conv\{z_{i}\}_{i \in I}$.

 On the simple examples, it will be further illustrated that the matrix \,$B$\,
 may be singular in some cases.
 In this events, the function \,$\varphi(\alpha)$\, is not strictly convex a fortiori strongly convex.
\begin{examp}\label{P_1}
 Let \,$z_{1}=(1,0), \,  z_{2}=(2,0)$,\, i.e. the vectors \,$z_{i}, i=\overline{1,2}$\, be linearly dependent. Then the matrix \,$B=$\,
 $\left(
\begin{array}{cc}
  1 & 2 \\
  2 & 4 \\
\end{array}
\right)$ is singular.
 \end{examp}
\begin{examp}\label{P_2}
 Let \,$z_{1}=(-1,0), \,  z_{2}=(1,0),$\, then \,${\bf 0} \in conv \bigl \{z_{i}\bigr \}_{i =\overline{1,2}}.$\,
 Having constructed the matrix \,$B = \left(
\begin{array}{cc}
  \enspace 1 & -1 \\
  -1  & \enspace 1
\end{array}
\right),$
 one can easily see that it is also singular.
\end{examp}
\begin{examp}\label{P_3}
Now, we consider the example for which it holds \,$m \gg n$.\, Let
\,$n = 1,$\, $m = 3,$\, \,$z_{1} = -1,$\, \,$ z_{2} = 1,$ $z_{3} =
5$.\, In this case,  \,$C = \left( \begin{array}{c}
  \enspace 1  \\
  -1  \\
  -5
\end{array}
\right).$\, Consequently, the~matrix \,$C \cdot C^{T} = \left(
\begin{array}{ccc}
   \enspace 1 & -1 & -5 \\
  -1 &  \enspace 1 &  \enspace 5 \\
  -5 &  \enspace 5 & 25
\end{array}
\right)$ \, \medskip \medskip is singular, because \,$det\,C \cdot
C^{T} = 0.$
\end{examp}
By construction, any point from \,$L:= conv \bigl \{z_{i}\bigr
\}_{i \in
  I}$\, is representable as a convex
combination of the vectors \,$z_{i}$.\, This implies that the
problem~(\ref{G_norm1})--(\ref{G_norm3})
   is equivalent to the program~(\ref{t_1}).

  Let \,$\alpha^{*}=(\alpha^{*}_{1}, \alpha^{*}_{2}, \ldots,
  \alpha^{*}_{m})$\, be an optimal solution of the
  problem~(\ref{G_norm1})--(\ref{G_norm3}). Then, the following
  statement is evident.
\begin{lemm}\label{A_13}
\quad $\rho = \sum \limits_{i \in I} \alpha^{*}_{i}z_{i}.$
\end{lemm}

We remind that \,$\rho ={\bf 0}$\, if and only if \,${\bf 0} \in
L$.\, As a consequence of this fact, Lemma~\ref{A_13}\, gives the
obvious justification of the next result.
\begin{lemm}\label{A_14}
\quad $\varphi(\alpha) = {\bf 0}$\,\, if and only if  \,${\bf 0}
\in L.$\,
\end{lemm}

On account of the specifity of the problem setting
for~(\ref{G_norm1})--(\ref{G_norm3}), we formulate an optimality
criterion of its solutions. Having gotten the solution $\alpha$,\,
we need to check the fulfillment of the inequalities
$$\langle z_{i}, \sum \limits_{j=1}^{m}
\alpha_{j}z_{j}\rangle \geq \|\sum \limits_{j=1}^{m}
\alpha_{j}z_{j}\|^{2}= \langle \alpha, B\alpha \rangle \quad
\forall i=\overline{1,\, m},\, \mbox{ т.е.} \sum \limits_{j=1}^{m}
\alpha_{j}b_{ij} \geq \langle \alpha, B\alpha \rangle \enspace
\forall i = \overline{1,\, m},$$ where \,$b_{ij}=\langle
z_{i},z_{j}\rangle.$\,

Thus, for solving the problem of projecting the origin onto the
convex polyhedron \,$L$,\, it suffices to find the solution
\,$\alpha^{*}$\,to the problem~(\ref{G_norm11})--(\ref{G_norm31}).
Then, according to Lemma~\ref{A_13}, we calculate \,$\rho=\sum
\limits_{i \in I} \alpha_{i}^{*}z_{i}.$\,

Let us carry out some comparison of the problems~(\ref{lagr1})\,
and~(\ref{G_norm1})--(\ref{G_norm3}). Both problems have the
constraint on the nonnegativity of their variables. The
problem~(\ref{lagr1}),\, in contrast to
(\ref{G_norm1})--(\ref{G_norm3}), has no constraints except one
which ensures that its decision variables are nonnegative.

The functions  \,$f(u)$\, and  \,$\varphi(\alpha)$\, have the same
squared  part \,$\langle x, C C^{T} x \rangle.$\, In the
problem~(\ref{G_norm1})--(\ref{G_norm3}), the linear function
\,$\langle e, x \rangle$\,  involved in the system of constraints,
while in the problem~(\ref{lagr1}) it is in the objective
function. The function \,$\varphi(\alpha)$\, is bounded from
below: \,$\varphi(\alpha)\geq 0 \quad \forall \alpha \geq 0.$\,
The equality \,$ \, \varphi(\alpha)= 0$\, is fulfilled if and only
if
 \,${\bf 0} \in L.$\,
The same condition is in fact necessary and sufficient for the
function \,$f(u)$\, to be unbounded from below on the domain of
the feasible solutions.

Finally, we note that in \cite{G_Stets_Nurmi}, there was modified
the problem (\ref{G_norm1})--(\ref{G_norm3}) with the purpose of
getting an equivalent nonsmooth problem and further effective
solving it by  the nonsmooth penalty  and subgradient algorithms.

\section{Reduction to the Maximin Problem}
%\section*{\centerline {$5^{0}. \quad \mbox{Седьмая постановка}$}}

This section is devoted to the investigation of how the reduction
of the learnt projection problem can be made to the maximin
problem. Solving of the maximin problem may be carried out by
simple transition to the minimax problem and the subsequent
application of the software package. For instance, a package
Optimization Toolbox in MATLAB  includes the function fminimax
which is used to solve the minimax constraint problem.

 Let \,$c^{*}$\, be a solution of the program
\begin{equation}\label{G_problem1}
    \max \limits_{\|c\| = 1} \,t_{L}(c),
\end{equation}
where \,$ t_{L}(c) = \inf \limits_{z \in L} \,\langle
 c,z\rangle.$\, By virtue of continuity on the whole space,
the linear function $\langle c,x\rangle$, $c \in \mathbb{R}^{n}$
attains its  infimum on the compact set $L$. Consequently, it
holds $\inf \limits_{z \in L} \,\langle
 c,z\rangle = \min \limits_{z \in L} \,\langle
 c,z\rangle.$\,  Let  $z^{*}$ be the optimizer of the following problem
                       \begin{equation}\label{PR_1}
                              \min \limits_{z \in L} \,\langle c,z\rangle
                       \end{equation}
                       for some $c
\in \mathbb{R}^{n}$. Note that, from the above, it obviously
follows \,$z^{*} \in L$.
\begin{lemm}\label{A_15_prim}
If  $z^{*} \in L$ is the solution of~(\ref{PR_1}), then $\langle
c,z^{*}\rangle = \min \limits_{i \in I} \, \langle
c,z_{i}\rangle.$
\end{lemm}
\proof Suppose that the assertion of the lemma is not true, i.e.
for all $i \in I$ it takes place the correlation $\langle
c,z_{i}\rangle
> \langle c,z^{*}\rangle$. By definition of the set  $L$, any
 its point can be performed as the convex combination of
 $z_{i}, \, i \in I$. Then, it is fulfilled
$$\langle c,z^{*}\rangle = \langle c, \sum \limits_{i \in I}
\alpha_{i}z_{i}\rangle =  \sum \limits_{i \in I} \alpha_{i}\langle
c,z_{i}\rangle > \langle c,z^{*}\rangle,\, \mbox{где} \enspace
\alpha_{i} \geq 0, \sum \limits_{i \in I} \alpha_{i} = 1.$$ The
obtained contradiction proves that our assumption is not true,
i.e. there can be found at least one index $i_{0} \in I$ such that
 $\langle c,z_{i_{0}}\rangle = \langle c,z^{*}\rangle$.
\,$\square$

 The statement of the preceding lemma means that a minimum of the linear function \,$\langle c,z\rangle,$\,$c \in \mathbb{R}^{n}$\,
 on the polyhedron, specified as the convex hull of a finite number
 of points, is reached at least at one of the points.
  Then, for the problem~(\ref{G_problem1}), on account of the next
  equality
$$\min \limits_{z \in L} \,\langle
 c,z\rangle = \min \limits_{i \in I} \,\langle c,z_{i}\rangle, \, \forall c \in  \mathbb{R}^{n},$$
we have \,$ t_{L}(c) = \min \limits_{i \in I} \langle
 c,z_{i}\rangle.$\,
\begin{lemm}\label{A_15}
1.\, If \,$t_{L}(c^{*}) > 0,$\, $v = c^{*} \cdot
t_{L}(c^{*}),$\, then it holds \,${\bf P}_{L}({\bf 0}) = v.$\\
 2.\, If \,$t_{L}(c^{*}) \leq 0,$\, then
 \,${\bf P}_{L}({\bf 0}) = {\bf 0}.$
\end{lemm}
   The proof of the previous lemma presented in \cite{G_Gabid6}
   (see the proof of Lemma 3.8, p.~167). In \cite{G_Gabid6}, this statement was proved for more general setting, namely for the case
    when the set, onto which there is carried out the projection, is convex and closed on
\,$\mathbb{R}^{n}.$ In \cite{G_Gabid6}, there is justified the
solvability  of the problem~(\ref{G_problem1}) (see p.149) and
uniqueness of its solution (see p.167). The
problem~(\ref{G_problem1}) may be reduced to the following
problem:
\begin{equation}\label{G_problem2}
    \max \limits_{\|c\| \leq 1} \,t_{L}(c).
\end{equation}
 The interconnection of the solutions for problems~(\ref{G_problem1}) and (\ref{G_problem2}) is described  in the next lemma.
\begin{lemm}\label{A_16}
If \,$\hat{c} \neq {\bf 0},$\, and  \,$c^{*}$\ are the solutions
of the problems~(\ref{G_problem2}) and  (\ref{G_problem1}),
respectively, then it holds $$t_{L}(c^{*}) =
t_{L}\left(\frac{\hat{c}}{\|\hat{c}\|}\right).$$
\end{lemm}
For details of the lemma's proof in more general setting (for the
case when
 $L $ is any convex and closed subset of $\mathbb{R}^{n}$), we refer the interested reader to \cite{G_Gabid6}, p.167. Thus, the
previous lemma proves that it is possible  to reduce the
problem~(\ref{G_problem1}) to~(\ref{G_problem2}).

For justifying of some additional results, we shall use the lemma
on the necessary and sufficient conditions for emptiness of the
cone of generalized strong support vectors  of  the set \,$L$:
\\\centerline{$V_{L}:= \{c \in \mathbb{R}^{n}: \min \limits_{z \in
L} \,\langle c,z\rangle
> 0\}$.}
Due to Lemma~\ref{A_15_prim}, it is obviously fulfilled \,$V_{L} =
\{c \in \mathbb{R}^{n}: \min \limits_{i \in I} \,\langle
c,z_{i}\rangle
> 0\}$.
          \begin{lemm}\label{A_17_prim}
 $V_{L} = \emptyset$\, if and only if \,${\bf 0} \in L$.
          \end{lemm}
          The statement of the lemma is the particular case for the statement of Lemma 3.7\,
          from \cite{G_Gabid7}.
\begin{lemm}\label{A_17}
If \,${\bf 0}\notin L,$\, \,$\hat{c}$\, is the solution of~(\ref{G_problem2}), then there are fulfilled the following statements:\\
1.\, \,$t_{L}(\hat{c}) > 0,$\\
2.\, the solution of the problem~(\ref{G_problem2}) is unique.
\end{lemm}
\proof 1.\,   Due to Lemma~\ref{A_17_prim},\,  since \,${\bf
0}\notin L$,\, we have  \,$V_{L} \neq \emptyset$,\, i.e. the
origin of $\mathbb{R}^{n}$
 is strongly separable from the polyhedron $L$. Theorem~2.1 from \cite{G_Gabid6}\,
asserts that for the strong separability between the origin and
some nonempty set $\Phi \subset \mathbb{R}^{n}$,\, a necessary and
sufficient condition is holding \,$t_{\Phi}(c^{*})
> 0,$ where  \,$c^{*}$ is the solution of the problem $\max \limits_{\|c\| = 1}
 \,t_{\Phi}(c),$\, $ t_{\Phi}(c) = \inf \limits_{z \in \Phi} \,\langle
 c,z\rangle.$\, Due to this theorem, we observe the implication  \,$t_{L}(c^{*}) > 0 \,
 \Rightarrow \,\max \limits_{\|c\| \leq 1}
 \,t_{L}(c) \geq \max \limits_{\|c\| = 1}
 \,t_{L}(c) > 0$.

2.\, Suppose that the solution of the problem~(\ref{G_problem2})
is not unique. Let for the vector  \,$c_{1}$\, be fulfilled
\,$\|c_{1}\| \leq 1$,\, $t_{L}(c_{1}) = t_{L}(\hat{c}),$\, $c_{1}
\neq \hat{c}$.\, By the fact that \,$t_{L}(\hat{c}) > 0$\,
together with  $t_{L}(c_{1}) = t_{L}(\hat{c}),$\, it holds
$\hat{c} \neq {\bf 0}, \,$ $c_{1} \neq {\bf 0}.$ Then from
Lemma~\ref{A_16}, we get \,$t_{L}\left(\frac{\textstyle
c_{1}}{\textstyle \|c_{1}\|}\right) = t_{L}\left(\frac{\textstyle
\hat{c}}{\textstyle \|\hat{c}\|}\right),$\, $\frac{\textstyle
c_{1}}{\textstyle \|c_{1}\|} \neq \frac{\textstyle
\hat{c}}{\textstyle \|\hat{c}\|}$.\, Clearly this contradicts the
uniqueness of the solution  for~(\ref{G_problem1}). The obtained
contradiction
 proves that our assumption is not true, i.e. the solution of~(\ref{G_problem2})
is really unique.\,$\square$

\begin{lemm}\label{A_18}
For \,${\bf 0}\in L$\, to be fulfilled, it is  necessary and
sufficient to have  \,$t_{L}(\hat{c}) = 0.$
\end{lemm}
\proof Sufficiency can be proved by reductio ad absurdum. Suppose
that \,${\bf 0} \notin L.$\, Then, by force of Lemma~\ref{A_17},\,
we have \,$t_{L}(\hat{c})
> 0.$\, This contradicts to the condition of the lemma. Consequently, our assumption is not
true, i.e. it holds \,${\bf 0} \in L.$\,

Necessity. Owing to Lemma~\ref{A_17_prim},\, the statement that
\,${\bf 0}\in L$\, is the same as saying that the cone of
generalized strong support vectors for \,$L$\, is empty. This, in
its own turn, means that the origin of $\mathbb{R}^{n}$ is not
strongly separable from the polyhedron  \,$L.$\, We apply further
the afore-mentioned Theorem~2.1 from \cite{G_Gabid6}\, (see the
proof of the first statement of the previous lemma). Owing this
theorem,\, we have \,$t_{L}(c^{*}) \leq 0.$ By the formulation the
problem~(\ref{G_problem2}), it always holds \,$t_{L}(\hat{c}) \geq
0.$ Supposing that \,$t_{L}(\hat{c}) > 0$, by Lemma~\ref{A_16},
one obtains the contradiction with the inequality \,$t_{L}(c^{*})
\leq 0.$ This contradiction allows to complete the proof.
\,$\square$

         \begin{lemm}\label{A_19}
               If  \,$t_{L}(\hat{c}) > 0,$\, then ${\bf 0}\notin L.$\,
         \end{lemm}
\proof Suppose that the assertion of the lemma is not true, i.e.
${\bf 0}\in L.$\, Then, due to Lemma~\ref{A_18},\,
 it holds \,$t_{L}(\hat{c}) = 0.$\, This contradicts the condition of the lemma.
\,$\square$
                  \begin{lemm}\label{A_191}
               Let the vector $\hat{c} \neq {\bf 0}$ be the solution of the problem~(\ref{G_problem2}), then
                $\hat{c}$ belongs to surface of the unit ball, i.e. $\|\hat{c}\| = 1$.
                        \end{lemm}
\proof Suppose the converse, i.e. $\|\hat{c}\| < 1$. According to
Lemma~\ref{A_16}, we then obtain
         \begin{equation}\label{Eq_A191}
    t_{L}(c^{*}) =
t_{L}\left(\frac{\hat{c}}{\|\hat{c}\|}\right) =
\frac{1}{\|\hat{c}\|} \cdot t_{L}(\hat{c}) > t_{L}(\hat{c}).
         \end{equation}
We note that the vector $c^{*}$ belongs to the domain of the
feasible solutions for~(\ref{G_problem2}). Therefore, the
inequality~(\ref{Eq_A191}) contradicts to  the fact that
\,$\hat{c}$\, is the solution of~(\ref{G_problem2}). The obtained
contradiction allows to see that our assumption is not true, i.e.
it holds $\|\hat{c}\| = 1$. \,$\square$

 From Lemmas~\ref{A_16},  \ref{A_17}, \ref{A_191},\, we get the following corollary.
             \begin{seqv}\label{A_1922}
Let the vector \,$\hat{c} \neq {\bf 0}$ be the solution
of~(\ref{G_problem2}), then \,$c^{*} = \hat{c}$.
             \end{seqv}

             Note that the previous Lemmas~\ref{A_17},\,\ref{A_19} describe
a necessary and sufficient condition for holding \,${\bf 0}\notin
L.$\, Lemmas~\ref{A_16}--\ref{A_19} provide the possibility of
using the problem~(\ref{G_problem2}) (instead the
problem~(\ref{G_problem1})) for projecting the origin onto the
convex polyhedron. The problem~(\ref{G_problem2}) compares
favorably with the problem~(\ref{G_problem1}) in the fact that,
unlike the problem~(\ref{G_problem1}), the objective function
in~(\ref{G_problem2}) is optimized on the whole $n-$dimensional
unit ball, not only on its surface. As is widely known, the
numerical methods usually
 better work in the case of the problems with the inequality constraints, not
 with the equality ones.

 Owing to Lemma~\ref{A_15}\, and Corollary~\ref{A_1922}, the assertion of the following lemma is true.

\begin{lemm}\label{A_192} 1.\, If \,$t_{L}(\hat{c}) > 0,$\, $v =
\hat{c} \cdot
t_{L}(\hat{c}),$\, then it holds \,${\bf P}_{L}({\bf 0}) = v.$\\
 2.\, If  \,$t_{L}(\hat{c}) = 0,$\, then
 \,${\bf P}_{L}({\bf 0}) = {\bf 0}.$
\end{lemm}

\section{Reduction to the Linear Complementarity Problem}
%\section*{\centerline {$5^{0}. \quad \mbox{Седьмая постановка}$}}

 In this section, we realize the reduction of the problem of projecting the origin  onto the convex
 polyhedron to the linear complementarity problem (LCP).

Using the $n \times n$--dimensional unit matrix $E$, we rewrite
first the quadratic programming
problem~(\ref{Cye_1})--(\ref{Cye_2}) in the vector-matrix form

           \begin{equation}\label{Cye1_matrix}
                  min \, \langle y, Ey \rangle,
           \end{equation}
         \begin{equation}\label{Cye2_matrix}
                      - Cy \geq e.
         \end{equation}
         We need to convert the problem~(\ref{Cye1_matrix})--(\ref{Cye2_matrix}) to the canonical form.
The following form of the quadratic programming problem is
considered to be canonical:
 \begin{equation}\label{Cye3_matrix}
                  min \, f(\hat{x}) =  \langle \hat{x}, \hat{D}\hat{x} \rangle + \langle \hat{p}, \hat{x} \rangle,
           \end{equation}
         \begin{equation}\label{Cye4_matrix}
                      \hat{A}\hat{x} \geq \hat{b},
         \end{equation}
         \begin{equation}\label{Cye5_matrix}
                      \hat{x} \geq 0,
         \end{equation}
where $\hat{D}$  is an $n \times n$--dimensional matrix,\,
$\hat{A}$ is an $l \times n$--matrix of the constraint
coefficients, $\hat{b}$ is the  $l$--dimensional vector,
$\hat{p},\, \hat{x} \in \mathbb{R}^{n}$.

         Carry out the change of the variables in the problem~(\ref{Cye1_matrix})--(\ref{Cye2_matrix}):
         $y_{1} = s_{1} - s_{1}^{\,\prime},$ $y_{2} = s_{2} - s_{2}^{\,\prime},$ $\ldots,$ $y_{i} = s_{i} - s_{i}^{\,\prime},$
          $\ldots,$ $y_{n} = s_{n} - s_{n}^{\,\prime}$\, using the nonnegative variables $s_{i},\, s_{i}^{\,\prime},\, i = 1, 2, \ldots, n.$
Then, it is easy to see that the objective function takes a form:
         $$\langle y, Ey \rangle = s_{1}^{2} - 2s_{1}s_{1}^{\,\prime} + (s_{1}^{\,\prime})^{2}+ \ldots + s_{n}^{2} - 2s_{n}s_{n}^{\,\prime} + (s_{n}^{\,\prime})^{2} =
                    \langle x, Dx \rangle + \langle p,x \rangle,$$
        where $x \in \mathbb{R}^{2n},$ $x = (s_{1}, s_{2}, \ldots, s_{n}, s_{1}^{\,\prime}, s_{2}^{\,\prime}, \ldots,
        s_{n}^{\,\prime}),$ $p = \underbrace{(0, \ldots, 0)}_{2n},$
         $D$ has the dimension $2n \times 2n,$
       \,$D = \left(
\begin{array}{cc}
   \enspace \,E & \mid -E  \\
   \hline
   \enspace -E & \mid \enspace \enspace E  \\
  \end{array}
\right).$ \, Clearly, the matrix $D$ is symmetrical, i.e. $D =
D^{T}.$ We further substitute the variables in the system of
constraints in the problem
(\ref{Cye1_matrix})--(\ref{Cye2_matrix}):\\ $- C \cdot y = - C
\cdot \left(
      \begin{array}{c}
  s_{1} - s_{1}^{\,\prime} \\
  s_{2} - s_{2}^{\,\prime} \\
  $\vdots$ \\
                    s_{n} - s_{n}^{\,\prime} \\
\end{array}
\right) =$ $A \cdot x \geq b,$\, where the matrices \,$C$\, and
$A$\, have the dimension  $m \times n$ and  $m \times 2n$,
respectively,  $A = \left(
\begin{array}{cc}
-C \mid C
\end{array}
\right)$, $b = \underbrace{(1, \ldots, 1)}_{m},$ $x \geq 0.$ Under
the change of variables,  the problem
(\ref{Cye1_matrix})--(\ref{Cye2_matrix}) is transformed to the
following canonical form:
\smallskip
\begin{equation}\label{Cye3_matrix1}
                  min \, f(x) =  \langle x, Dx \rangle + \langle p, x \rangle,
           \end{equation}
         \begin{equation}\label{Cye4_matrix1}
                      Ax \geq b,
         \end{equation}
         \begin{equation}\label{Cye5_matrix1}
                      x \geq 0.
         \end{equation}
\begin{examp}
  For $n = 2$, the matrix of the quadratic form $\langle x,Dx \rangle$
can be written as: \\ \smallskip \centerline{$D = \left(
\begin{array}{cccc}
   \enspace 1 & \enspace 0  &\,\mid -1 & \enspace 0    \\
  \enspace 0 & \enspace 1   &\mid \enspace 0 & -1   \\
     \hline
    \enspace -1 & \enspace 0 &\mid \enspace 1 & \enspace 0   \\
  \enspace 0 & \enspace  -1 &\mid \enspace 0 & \enspace 1   \\
    \end{array}
\right).$}
\end{examp}
 As is known, the type of the problem~(\ref{Cye3_matrix1})--(\ref{Cye5_matrix1})
 depends on the sign-definiteness of the matrix corresponding to the quadratic
 form $\langle x, Dx \rangle$. If the matrix $D$ is non-negatively (non-positively)
definite, then (\ref{Cye3_matrix1})--(\ref{Cye5_matrix1})\, is the
convex (nonconvex) programming problem. In the case of the
positive-definiteness of the matrix $D$, the objective function of
the problem is strongly convex.
 In this event, the problem
(\ref{Cye3_matrix1})--(\ref{Cye5_matrix1}) has the unique
solution.

 Let us remind (see \cite{G_Gantmakher}) that the real-valued quadratic form $\langle x,
\hat{D}x \rangle$ is called positively definite, if for real $x
\neq {\bf 0}$ it holds $\langle x, \hat{D}x \rangle
> 0$.
The real-valued quadratic form is called non-negatively definite
if for $x \in \mathbb{R}^{n}$ it is fulfilled $\langle x, \hat{D}x
\rangle \geq 0$ . The positive (nonnegative) definiteness of the
quadratic form corresponds to the same sign-definiteness of the
associated  matrix.

We further clarify what is the sign-definiteness of the block-structural matrix \\
\centerline{$D = \left(
\begin{array}{cc}
   \enspace \,E & \mid -E  \\
   \hline
   \enspace -E & \mid \enspace \enspace E  \\
  \end{array}
\right).$} \\
By construction, $\langle x, Dx \rangle = (s_{1} -
s_{1}^{\,\prime})^{2} + \ldots + (s_{n} -s_{n}^{\,\prime})^{2}
\geq 0.$  For all cases, even when it simultaneously holds $s_{i}
\neq 0,\, s_{i}^{\,\prime} \neq 0, i =1, \ldots, n,$ but $s_{i} =
s_{i}^{\,\prime}, \forall i =1, \ldots, n$ we will obtain $\langle
x, Dx \rangle = 0$. Therefore, $\langle x, Dx \rangle$ satisfies
 the definition of the non-negatively definite quadratic form. As
 a consequence, the matrix $D$ is also nonnegative definite.

Write the optimality conditions of Kung--Tucker for the problem of
minimizing the convex function subject to the linear constraints
((\ref{Cye3_matrix1})--(\ref{Cye5_matrix1})): it is necessary to
find the vectors $x,\, z,\, u,\, r$\, such that
$$u = 2Dx - A^{T}z + p^{T} = 2Dx - A^{T}z,$$
$$r = Ax -b,$$
$$x,\, z,\, u,\, r > 0,$$
$$u^{T}x + r^{T}z = 0.$$
These are the necessary and sufficient conditions for the solution
of the problem (\ref{Cye3_matrix1})--(\ref{Cye5_matrix1}) to be
optimal. Following to \cite{G_Cottle_Dantzig}, we represent this
system as LCP which consists in finding the vectors $w$\, and
$v$\, such that

\begin{equation}\label{G_dop1}
  w = Mv + q,
\end{equation}
{\begin{equation}\label{G_dop2} w \geq 0, v \geq 0,
\end{equation}
\begin{equation}\label{G_dop3}
 \langle w , v \rangle = 0,
\end{equation}}
 where $M$\, is  an $k \times k$--dimensional matrix;\enspace
$w, v, q$\,  are the  $k$--dimensional vectors, \smallskip \\
 \centerline{$w =
\left(
\begin{array}{c}
   u \\
      r \\
  \end{array}
\right),$ $v = \left(
\begin{array}{c}
   x \\
      z \\
  \end{array}
\right),$} $M = \left(
\begin{array}{cc}
   \enspace \,2D & \, \mid  -A^{T}  \\
   \hline
   \enspace A & \mid \enspace \enspace \Theta  \\
  \end{array}
\right),$ $k = 2n + m,$ $ q = \left(
                 \begin{array}{c}
    \, 0 \\  \vdots \\  \, 0 \\ -  b
                  \end{array}
  \right),$ $\Theta$\, is the null matrix with dimension $m \times
  m$.\smallskip

  It is known that to any optimal solution of LCP (by virtue of the connection through the conditions of Kung--Tucker)
  corresponds the optimal solution of the quadratic programming problem. However, the converse is not true (except the case when the value of the
  objective function of the quadratic programming problem is equal to zero).
   In \cite{G_Cottle}, there are studied  the conditions  to which
   must satisfy
   the matrix $M$ in order to guaranteed that the solutions of LCP
   are identical to the points of Kung--Tucker for
the quadratic programming problem (associated with this LCP). It
is known that if the matrix
 $D$ is non-negatively definite, then  the matrix $M$ is also
non-negatively definite \cite{G_Cottle_Dantzig}. For this case, it
is proved that if there exists the solution of the complementarity
problem, then  Lemke's method \cite{G_Lemke} allows to find it for
the finite number of iterations.  If the linear complementarity
problem is without solutions, then the quadratic programming
 problem has no solutions, too.  In our case, this means that the
 domain of the feasible solutions for the problem
(\ref{Cye1_matrix})--(\ref{Cye2_matrix}) is empty. According to
Lemma~\ref{l9}, we then have \,${\bf 0} \in L$.\,

 The results of the experimental research published in
\cite{G_Ravindran}\, have shown that, for the case when the matrix
$\hat{D}$ of the problem~(\ref{Cye3_matrix})--(\ref{Cye5_matrix})
 is non-negatively definite,  Lemke's method has  several advantages in comparison to
the majority of the  quadratic programming methods. Let us note
that the function LCPSolve in MATLAB solves the linear
complementarity problem using a pivoting algorithm.

Further, we briefly consider the representation of the
optimization problem (\ref{G_norm11})--(\ref{G_norm31}) as LCP
(\ref{G_dop1})--(\ref{G_dop3}). In this case, $M = \left(
\begin{array}{cc}
   \enspace \,2B & \, \mid -A^{T}  \\
   \hline
   \enspace A & \mid \enspace \enspace \Theta  \\
  \end{array}
\right),$ $ q = \left(
                 \begin{array}{c}
     0 \\ \vdots \\  0 \\  -1 \\ 1
                  \end{array}
  \right),$
  where $M$ is a~quadratic matrix with dimension $(m+2) \times (m+2)$;\enspace
$w, v, q$ are the $(m + 2)$--dimensional vectors, $\Theta$ is the
null matrix of dimension $2 \times 2$,\, $B = C \cdot C^{T},$\, $A
= \left(
                 \begin{array}{c}
     e^{T} \\ -e^{T}
                  \end{array}
  \right).$

The problem~(\ref{lagr1})  may also be represented  as
LCP~(\ref{G_dop1})--(\ref{G_dop3})
 with  $$M =\frac{\textstyle 1}{\textstyle 2}B,
 \enspace
  q = -e,$$
   $M$\, is an $m \times m$--dimensional matrix;\,
$w, v, q$\,  are the \,$m$--dimensional vectors.

Thus, for the considered three quadratic programming problems, LCP
associated with the problem~(\ref{lagr1}) compares favorably with
the others. The advantage of this complementarity problem  is that
it has the least dimension.

\section{Reduction to the Nonnegative Least Squares
Problem.}
%\section*{\centerline {$5^{0}. \quad \mbox{Восьмая постановка}$}}

  In this section, we investigate  how the problem of projecting the origin onto the convex polyhedron
 can be reduced to the nonnegative least squares
problem (NLSP).

 In that scientific and practical applications, where one needs to estimate some
 vector of observations
  \,$b \in \mathbb{R}^{m}$\, by the \,$n$\,  basis
 factors or measures \\ \centerline{$A_{i}, \, i = \overline{1,n},$\, $A=[A_{1},A_{2}, \ldots, A_{n}] \in \mathbb{R}^{m\times
 n},$}
 there arises  the least squares
problem. The nonnegative least squares problem
 is usually formulated as the following optimization problem:
 \begin{equation}\label{lsq1}
   \min_{x\geq 0}\, \frac{1}{4} \|Ax-b\|^{2},
\end{equation}
$$ A \in \mathbb{R}^{m\times n},\, b \in \mathbb{R}^{m},$$
 which consists in finding the nonnegative vector
$x$ maximizing the closeness of the vector of deviations to the
origin of the space (in sense of the Euclidean distance).
  Clearly, the problem~(\ref{lsq1}) may be converted to the quadratic programming problem with the nonnegative variables:
\begin{equation}\label{lsq2}
   \min_{x\geq 0}\,\, f(x) = \frac{1}{4} x^{T}Qx + p^{T}x,
\end{equation}
$$\mbox{where}\enspace Q = A^{T}A, p = -\frac{1}{2}A^{T}b.$$
Let us note that (\ref{lsq1}) is the convex programming
 problem, since the matrix \,$Q = A^{T}A$\,  is non-negatively definite and the non-negativity constraints
specify the convex set of the feasible solutions $x$.

To establish a connection of the problem of projecting the origin
onto the convex polyhedron with NLSP, we further consider the
problem~(\ref{lagr1}). The problem~(\ref{lagr1}) may be reduced
to~(\ref{lsq2})
 under conditions:
\begin{equation}\label{cond12}
1)\, A = C^{T},\quad 2)\, A^{T}b =2e.
\end{equation}

 On the simple examples, we further demonstrate  that it is not hard to achieve
  the fulfillment of the second condition
 in~(\ref{cond12}).
 \begin{examp}\label{exp1}
 $\enspace A = \left(
  \begin{array}{cc}
   \enspace 1 & \enspace -1\\
  \enspace 0 & \enspace \,1\\
    \end{array}
    \right),$
     $ A^{T} = \left(
  \begin{array}{cc}
   \enspace \,1 & \enspace 0\\
  \enspace -1 & \enspace 1\\
    \end{array}
    \right),$
     $ b = \left(
  \begin{array}{c}
   \enspace 2 \\
  \enspace 4 \\
    \end{array}
    \right)$.
        \end{examp}

\begin{examp}\label{exp2}
 $\enspace A = A^{T} = \left(
  \begin{array}{cc}
   \enspace -10 & \enspace 0\\
  \enspace \,\,\,0 & \enspace 5\\
    \end{array}
    \right),$
         $ b = \left(
  \begin{array}{c}
   \enspace  -0.2 \\
  \enspace  \,\,\,0.4 \\
    \end{array}
    \right)$.
        \end{examp}

Below, we consider the rule of construction of the vector $b$
satisfying the second condition  of~(\ref{cond12}) for the case of
the quadratic diagonal matrix  $A \in \mathbb{R}^{n\times n}$
(with zero non-diagonal elements and nonzero elements of the
principal diagonal):
$$b_{i} = \frac{2}{a_{ii}}, \, i = \overline{1,n},\, a_{ii}\neq 0.$$
For this case, we consider the vector $x^{*}$\, having
 coordinates:
$$x_{i}^{*} = \frac{b_{i}}{a_{ii}} = \frac{2}{a_{ii}^{2}}, \, i = \overline{1,n},\, a_{ii}\neq 0.$$
By construction, we have $x_{i}^{*} > 0, \, i = \overline{1,n}.$
Consequently, the vector  $x^{*}$\, is the solution of the
equi\-va\-lent problems~(\ref{lsq1})\, and (\ref{lsq2}). If the
vector $b$\, and matrix $A$ of the problem~(\ref{lsq2}) satisfy
the conditions~(\ref{cond12}), then the obtained vector $x^{*}$\,
is also the solution of the problem~(\ref{lagr1}). Then,\, the
solution of~(\ref{Cye_1})--(\ref{Cye_2})  may be calculated by the
following formula:
  \begin{equation}\label{y*}
y^{*}= -\frac{1}{2} C^{T}x^{*}.
\end{equation}
 We further determine the solution of the original
problem~(\ref{t_1}) by direct calculation as follows:
\,$\rho=y^{*}/\|y^{*}\|^{2}.$\,

Let $A = C^{T}$ be non-singular quadratic matrix of dimension
 $n\times n$. In this case, the vector  $b$\, satisfying the second condition of~(\ref{cond12}) can be determined as
 follows:\\
\centerline{$b = 2(A^{T})^{-1}e = 2C^{-1}e.$}
 Based on linear algebra, we may state that a necessary and sufficient condition for the matrix $A$\,
  to be nonsingular is that the given polyhedron $L$\,
   be specified
  by  the system of the linearly independent
   vectors $z_{i},\, i = \overline{1,n}.$

The first known method for solving NLSP has been proposed in
  \cite{G_Lowson_Hanson}. Nowadays, this method is known as the active-set method.
For instance, the afore-mentioned  method is realized in MATLAB by
help of function lsqnoneg. In the  literature,  there are  known
the fast modifications  and parallel variants of  the active-set
method \cite{G_Bro_Jong},  \cite{G_Luo_Duraiswami}. Comparing
analysis of the methods for solving NLSP can be found, for
instance, in
  \cite{G_Chen_Plemmons}.

%\section*{\centerline {$4^{0}. \quad \mbox{Пятая и шестая постановки}$}}
\section{Conclusions}
The main contributions of the present paper may be briefly
summarized as follows:\\
We have treated the opportunity to reduce PPOCP to the different
relevant problems of mathematical
 programming (MP) such as QPP, maximin problem, LCP, and NLSP.
Such reduction makes it possible to utilize a much more broad
spectrum of powerful tools of MP for solving PPOCP.

\end{document}